\documentclass[leqno]{amsart}

\usepackage{amsmath,amsfonts,amsthm,amssymb, mathrsfs, mathabx, accents, mathdots, textcomp, fontenc,  vntex}

\usepackage[english]{babel}

\usepackage[bookmarks]{hyperref}
\hypersetup{backref, colorlinks=true}

 \usepackage{graphicx}
\usepackage[all]{xy}
\newtheorem{theorem}{Theorem}[section]
\newtheorem{lemma}[theorem]{Lemma}
\newtheorem{proposition}[theorem]{Proposition}
\newtheorem{corollary}[theorem]{Corollary}
\newtheorem{remark}[theorem]{Remark}

\newtheorem{example}[theorem]{Example}
\newtheorem{definition}[theorem]{Definition}

\newtheorem{conjecture}[theorem]{Conjecture}
\newenvironment{equationth}{\stepcounter{theorem}\begin{equation}}{\end{equation}}
\newenvironment{preuve}{{\em{\noindent \textbf{Proof.} }}}
{\hfill $\blacksquare$}

\def\C{ \mathbb{C}}

\def\R{ \mathbb{R}}
\def\P{ \mathbb{P}}
\def\N{ \mathbb{N}}

\def\S{ \mathbb{S}}

\def\graph{ {\rm graph}}

\def\codim{ {\rm codim}}
\def\Rank{ {\rm rank}}
\def\corank{ {\rm corank}}
\def\Reg{ {\rm Reg}}
\def\Sing{ {\rm Sing}}
\def\rang{ {\rm Rank}}
\def\Rang{ {\rm Rank}}

\def\rond{\mathaccent"7017}

\begin{document}
\large


\title[]{ A REMARK ON THE TOPOLOGY OF A POLYNOMIAL MAPPING $F: \C^n \to \C^{n}$ VIA INTERSECTION HOMOLOGY}

\makeatother

\author[Nguy\~{\^e}n Th\d{i} B\'ich Th\h{u}y]{Nguy\~{\^e}n Th\d{i} B\'ich Th\h{u}y}
\address[{Nguy\~{\^e}n Th\d{i} B\'ich Th\h{u}y}]{UNESP, Universidade Estadual Paulista, ``J\'ulio de Mesquita Filho'', S\~ao Jos\'e do Rio Preto, Brasil}
\email{bichthuy@ibilce.unesp.br}
\maketitle \thispagestyle{empty}
\begin{abstract}
In  \cite{Valette}, Guillaume and Anna Valette associate singular varieties $V_F$ to a polynomial mapping $F: \C^n \to \C^n$. In the case 
%Singular varieties  $V_F $  contructed in \cite{Valette},  associated to a polynomial mapping $F : \C^{n} \to \C^{n}$, satisfy the following property:  in the case 
$F: \C^2 \to \C^2$, if the set $K_0(F)$ of critical values  of $F$ is empty, then  $F$ is not proper if and only if the 2-dimensional homology or  intersection homology (with any perversity)  of $V_F$ are not trivial. 
In \cite{ThuyValette}, the results of \cite{Valette} are generalized in the case $F: \C^n \to \C^n$ where $n \geq 3$, with an additional condition. In this paper, we prove that 
 if  $F: \C^2 \to \C^2$  is a non-proper {\it generic dominant } polynomial mapping,  then the 2-dimensional homology and  intersection homology (with any perversity)  of $V_F$ are not trivial. We prove that this result is true also for  a non-proper {\it generic dominant } polynomial mapping $F: \C^n \to \C^n$ ($\, n \geq 3$), with the same additional condition than in \cite{ThuyValette}. 
%then the results of \cite{Valette} and \cite{ThuyValette} hold also in  the case $K_0(F)$ is non empty. 
In order to compute the  intersection homology of the variety $V_F$, 
 we provide 
an explicit Thom-Mather stratification of the set $K_0(F) \cup S_F$. 
\end{abstract}
%\footnote{Ici, le mot "provide" est pr\'ecis? En fait, je ne "provide" pas, je montre qu'il existe une stratification de Thom-Mather compatible avec les stratifications de Thom-Mather de $S_F$ e $K_0(F)$ qui existent d\'ej\`a. Mais, de toute fa{\c c}on, le mot "provide" est aceptable pour le r\'esum\'e? Puisque, dans l'introduction, je dis plus clair, exactement ce que je fais.} 
\section{Introduction}

%In \cite{Valette}, the variety $V_F$ is constructed as follows: Let $F: \C^n \to \C^n$ be a polynomial mapping. We consider $F$ as a real one $F: \R^{2n} \to \R^{2n}$. By $Sing F$ we mean the set of critical points of $F$. 
%Thanks to Lemma 2.1 of \cite{Valette}, 
% there exists a covering $\{ U_1, \ldots , U_p \}$ of $M_F = \R^{2n} \setminus Sing(F)$ by open semi-algebraic subsets (in $\R^{2n}$) such that on every element of this covering, the mapping $F$ induces a diffeomorphism onto its image. We may find some semi-algebraic closed  subsets $V_i \subset U_i$ (in $M_F$) which cover $M_F$ as well. By the Mostowski's Separation Lemma (see \cite{Mos}, page 246), for each $ i =1, \ldots , p$, there exists a Nash function $\psi_i : M_F \to \R$,  such that  $\psi_i$ is positive on $V_i$ and negative on $M_F \setminus U_i$. 
% We can construct the Nash functions $\psi_i$ such that $\psi_i (x_k)$ tends to  zero when $\{x_k\} \subset M_F$ tends to infinity. 
%We define $V_F : = \overline{(F, \psi_1, \ldots, \psi_p)(M_F)},$ 
%that means, $V_F$ is the closure of the image of $M_F$ by $(F, \psi_1, \ldots, \psi_p)$. Then the variety $V_F$ is a real algebraic singular  variety of dimension $2n$, the singular set of which is contained in  $( S_F \times K_0(F)) \times \{0_\R\}^p$, where $S_F$ is the asymptotic set and $K_0(F)$ is the set of critical values of $F$ (proposition 2.3 of \cite{Valette}).  With this construction of $V_F$, we have the following property: in the case $F: \C^2 \to \C^2$ such that $K_0(F) = \emptyset$, 

In \cite{Valette}, Guillaume and Anna Valette provide a criteria for properness of a polynomial mapping $F: \C^2 \to \C^2$. They construct a real algebraic singular variety $V_F$ satisfying the following property: if the set of critical values of $F$ is empty then $F$ is not proper if and only if the 2-dimensional homology or  intersection homology (with any perversity)  of $V_F$ is not trivial (\cite{Valette}, theorem 3.2). 
In \cite{ThuyValette}, the result of \cite{Valette} is generalized in the general case $F: \C^n \to \C^n$, where $n \geq 3$, with an additional condition (\cite{ThuyValette}, theorem 4.5). The variety $V_F$ is a real algebraic  singular variety of dimension $2n$ in some $\R^{2n + p}$, the singular set of which is contained in  $(  K_0(F) \cup S_F ) \times \{0_\R\}^p$, where $K_0(F)$ is the set of critical values  and $S_F$ is the asymptotic set of $F$.

This paper proves that 
 if  $F: \C^2 \to \C^2$  is a non-proper {\it generic dominant } polynomial mapping,  then the 2-dimensional homology and  intersection homology (with any perversity)  of $V_F$ are not trivial. We prove that this result is true also for  a non-proper {\it generic dominant } polynomial mapping $F: \C^n \to \C^n$ ($\, n \geq 3$), with the same additional condition than in \cite{ThuyValette}. 
%if  $F: \C^n \to \C^n$ ($\, n \geq 2$) is a {\it generic dominant } polynomial mapping,  
%then the results of \cite{Valette} and \cite{ThuyValette} hold also in  the case $K_0(F)$ is non empty. 
To prove these results, we use the transversality theorem: if $F$ is non-proper {\it generic} dominant polynomial mapping, 
%\footnote{C'est la peine de donner la r\'ef\'erence ? Ou c'est classique et on n'a pas besoin de donner la r\'ef\'erence?}, 
we can construct an adapted $(2, \overline{p})$-allowable chain (in generic position)  providing  non triviality of  homology and  intersection homology of the variety $V_F$, for any perversity $\overline{p}$ (theorems \ref{Thuy1} and \ref{Thuy2}).  

In order to compute the intersection homology of the variety $V_F$ in the case $K_0(F) \neq \emptyset$,   we have to stratify the set $K_0(F) \cup S_F$. Furthermore, the intersection homology of the variety $V_F$ does not depend on the stratification  if we use a locally topologically trivial stratification. We know that a Thom-Mather stratification is a locally topologically trivial stratification (see \cite{Thom}, \cite{Whitney}). In \cite{Thom}, Thom defined  a partition of the set $K_0(F)$ by ``{\it constant rank}'', which is  a {\it local $\,$} Thom-Mather stratification. 
In \cite{Thuy}, the author provides a Thom-Mather stratification of the asymptotic set $S_F$ by  the {\it ``m\'ethode des fa{\c c}ons''}. 
In order to prove the main result, one important point of this paper is the following:  we show that in general the set $K_0(F)$ is not closed, so we cannot define a (global) stratification of $K_0(F)$ satisfying the frontier condition. Hence, we cannot define a (global) Thom-Mather stratification of $K_0(F)$. 
%\footnote{Est-ce qu'il faut dire pourquoi on a besoin the frontier conditions ? (pour une stratification de Thom-Mather. Ou, c'est \'evident comme cela ?}. 
However, we prove that the set $K_0(F) \cup S_F$ is closed and $K_0(F) \cup S_F = \overline{K_0(F)} \cup S_F$. This fact allows us to show that there exists a Thom-Mather stratification of the set $K_0(F) \cup S_F$ compatible with the partition of the set $K_0(F)$ defined by Thom and compatible with the stratification of the set $S_F$ defined by   the {\it ``m\'ethode des fa{\c c}ons''} (theorem \ref{prostraK0FSF}).

%\footnote{Pourquoi qu'il faut une stratification satisfait les conditions fronti\`eres pour $K_0(F)$?  Ok, $K_0(F)$ n'est pas ferm\'e mais m\^eme s'il est ouvert, par la definition d'une stratification, on peut toujours le stratifier, mais seulement que cette stratification n'est pas satisfait des conditions fronti\`eres, n'est-ce pas? Comment Thom a d\'efinie la stratification de 
%$K_0(F)$? 
%Est-ce que c'est une partition (stratification) que n'est pas satisfait des conditions fronti\`eres
% mais elle est une partition (stratification) de Thom-Mather? A discuter !!!}

%We provide also a Thom-Mather stratification of the set $K_0(F) \cup S_F$. 

%This paper provide also a Thom-Mather stratification of the set $K_0(F) \cup S_F$, using the Thom-Mather stratification of $K_0(F)$ defined by Thom (see, for example \cite{Thom}) and the Thom-Mather stratification of the asymptotic set $S_F$ defined by {\it ``la m\'ethode des fa{\c c}ons''} in \cite{Thuy}. We give also an example to illustrate the results of the paper. 

\section{PRELIMINARIES}
In this section we set-up  our framework. All the varieties we consider in this article are semi-algebraic.

%\subsection{Notations and conventions.}\label{section_notations}
%Given a topological space $X$, singular simplices of $X$ will be semi-algebraic continuous mappings $\sigma:T_i \to X$, 
%where $T_i$ is the standard $i$-simplex in $\R^{i+1}$. 
%Given a subset $X$ of $\R^n$ we denote by $C_i(X)$ the group of
%$i$-dimensional singular chains (linear combinations of singular simplices
%with coefficients in $\R$); if
%$c$ is an element of $C_i(X)$, we denote by $|c|$ its support.
%By $Reg(X)$ and $Sing(X)$ we denote respectively the regular and singular locus of the set $X$. Given $X \subset \R^n$, $\overline{X}$ will stand for the topological closure of $X$. The smoothness to be considered as the differentiable smoothness.
%
%Notice that, when we refer to the homology of a variety, the notation $H_*^{c}(X)$ refer to the homology with compact supports, the notation $H_*^ {cl}(X)$  refer to the homology with closed supports (see \cite{Jean Paul}).

\subsection{Intersection homology.}
We briefly recall the definition of intersection homology.  For details, 
we refer  to the fundamental work of M. Goresky and R. MacPherson
\cite{GM1} (see also \cite{Jean Paul}).

\begin{definition} 
{\rm Let $V$ be a $m$-dimensional semi-algebraic set.  A {\it semi-algebraic stratification of $V$} is the data of a finite semi-algebraic filtration 
$$V = V_{m} \supset V_{m-1} \supset \cdots \supset V_0 \supset V_{-1} = \emptyset,$$
 such that  for every $i$,  the set $V_i\setminus V_{i-1}$ is either an emptyset or a manifold of dimension $i$. A connected component of $V_i\setminus V_{i-1}$ is  called   {\it a stratum} of $V$.

Let $S_i$ be a stratum of $V$ and $\overline{S_i}$ its closure in $V$. 
If $\overline{S_i} \setminus S_i$ is the union of strata of $V$, for all strata $S_i$ of $V$, then we say that the stratification of $V$ satisfies the frontier condition. 
%\footnote{Why the frontier conditions are important?}
}
\end{definition}
%Notice that by  $\overline{S_i}$ we mean the closure of $S_i$ in $V$. 

\begin{definition} [see \cite{Thom}, \cite{Mather1}] \label{definitionthommather}
{\rm Let $V$ be a variety in a smooth variety $M$. 
 We say that a stratification of $V$ is a  {\it Thom-Mather stratification} 
 if each stratum  $S_i$ 
 is a differentiable variety of class $\mathcal C^{\infty}$  
 and if for each $S_i$, we have:

$a) \,$ an open neighbourhood (tubular neighbourhood) $T_i$ of $S_i$ in $M$,

$b) \,$ a  continuous retraction $\pi_i$ of $T_i$ on $S_i$,

$c) \,$ a continuous function  $\rho_i: T_i \to [0, \infty[$ which is  $\mathcal C^\infty$ 
 on the smooth part of $V \cap T_i$,

\noindent such that $S_i = \{ x \in T_i  : \rho(x) = 0 \}$ and if $S_i \subset \overline {S_j}$, then

\begin{enumerate}
\item[i)] the restricted mapping $(\pi_i, \rho_i): T_i \cap S_j \to S_i \times [0, \infty[$ is a smooth immersion,

\item[ii)] for $x \in T_i \cap T_j$ such that $\pi_j(x) \in T_i$, we have the following  relations of commutation:

$\qquad$ 1) $\pi_i \circ \pi_j(x) = \pi_i(x),$ 

$\qquad$ 2) $\rho_i \circ \pi_j(x) = \rho_i(x),$ 
\end{enumerate}

\noindent when the two members of theses formulas are defined. 
}
 \end{definition}

A Thom-Mather stratification of a variety $V$ satisfies the frontier condition. 

\medskip

We denote by $cL$  the open cone on the space $L$, the cone on the empty set being a point. Observe that if $L$ is a stratified set then $cL$ is stratified by  the cones over the strata of $L$ and an additional $0$-dimensional stratum (the vertex of the cone). 

\begin{definition}
{\rm A stratification of $V$ is said to be {\it locally topologically trivial} if for every $x \in V_i\setminus V_{i-1}$, $i \ge 0$, there is an open neighborhood $U_x$  of $x$ in $V$, a stratified set $L$ and a semi-algebraic homeomorphism 
 $$h:U_x \to (0;1)^i \times  cL,$$   such  that
 $h$ maps the strata of $U_x$ (induced stratification) onto the strata of  $  (0;1)^i \times cL$  (product stratification).}
\end{definition}

\begin{theorem} [see \cite{Thom}, \cite{Whitney}] \label{doesnotdepend}
A  Thom-Mather stratification is a  locally topologically trivial stratification. 
\end{theorem}

The definition of perversities has originally been given by Goresky and MacPherson:
\begin{definition}
{\rm 
A {\it perversity} is an $(m+1)$-uple  of integers $\bar p = (p_0, p_1, p_2,
p_3,\dots , p_m)$ such that $p_0 = p_1 = p_2 = 0$ and $p_{\alpha+1}\in\{p_\alpha, p_\alpha + 1\}$, for $\alpha \geq 2$.

Traditionally we denote the zero perversity by
$\overline{0}=(0, 0, \dots,0)$, the maximal perversity by
$\overline{t}=(0, 0, 0,1,\dots,m-2)$, and the middle perversities by
$\overline{m}=(0, 0, 0,0,1,1,\dots, [\frac{m-2}{2}])$ (lower middle) and
$\overline{n}=(0, 0, 0,1,1,2,2,\dots ,[\frac{m-1}{2}])$ (upper middle). We say that the
perversities $\overline{p}$ and $\overline{q}$ are {\it complementary} if $\overline{p}+\overline{q}=\overline{t}$.

Let $V$ be a semi-algebraic variety such that $V$ admits a  locally topologically trivial stratification. We say that a semi-algebraic subset $Y\subset V$ is  $(\bar
p, i)$-{\it allowable} if  
\begin{equationth} \label{conditionperversity}
\dim (Y \cap V_{m-\alpha}) \leq i - \alpha + p_\alpha \text{ for
all } \alpha \geq 2.
\end{equationth}
Define $IC_i ^{\overline{p}}(V)$ to be the $\R$-vector subspace of $C_i(V)$
consisting in the chains $\xi$ such that $|\xi|$ is
$(\overline{p}, i)$-allowable and $|\partial \xi|$ is
$(\overline{p}, i - 1)$-allowable.}
\end{definition}

\begin{definition} 
{\rm The {\it $i^{th}$ intersection homology group with perversity $\overline{p}$}, with real coefficients, denoted by
$IH_i ^{\overline{p}}(V)$, is the $i^{th}$ homology group of the
chain complex $IC^{\overline{p}}_*(V).$ 
}
\end{definition}

Notice that, the notation $IH_*^{\overline{p}, c}(V)$ refers to the intersection homology with compact supports, the notation $IH_*^ {\overline{p},cl}(V)$  refers to the intersection homology with closed supports. In the compact case, they coincide.
% and will be denoted by $IH_*^{\overline{p}}(V)$. 

 Goresky and MacPherson proved that the intersection homology  is independent on the
choice of the stratification satisfying the locally topologically trivial conditions \cite{GM1, GM2}.

The Poincar\'e duality holds for the intersection homology of a (singular) variety:

 \begin{theorem}[Goresky,
MacPherson \cite{GM1}]  \label{ThGM}
{\it For any orientable compact stratified
semi-algebraic $m$-dimensional variety  $V$,  the generalized Poincar\'e duality holds:
$$IH_k ^{\overline{p}}(V)  \simeq IH_{m-k} ^{\overline{q}} (V),$$
where $\overline{p}$ and $\overline{q}$ are complementary perversities.

 For the non-compact case, we have:
$$IH_k ^{\overline{p}, c}(V)  \simeq IH_{m-k} ^{\overline{q}, cl} (V).$$
}
\end{theorem}

\noindent

\subsection{The asymptotic set} \label{ensembleJelonek}

Let $F: \C^n \to \C^n$ be a polynomial mapping. Let us denote by $S_F$ the set of points at which  $F$ is non proper, {\it i.e.}, 
$$S_F: = \{ a \in \C^n \text{ such that } \exists \{ x_k\}_{k \in \N} \subset \C^n, \vert x_k \vert  \text{ tends to infinity and } F(x_k) \text{ tends to } a\},$$
where $ \vert x_k \vert$ is the Euclidean norm of  $x_k$ in $\C^n$. 
The set $S_F$ is called the asymptotic set of $F$. 

%Recall that, it is enough to define $S_F$  by considering sequences $\{\xi_k\}$ tending to infinity in the following sense: 
%each coordinate of these sequences either tends to infinity or converges.

In this paper, we will use the following important theorem:

\begin{theorem} \cite{Jelonek} \label{theoremjelonek1}
Let $F: \C^n \rightarrow \C^n$ be a polynomial mapping. 
 If $F$ is dominant,  {\it i.e.},  $\overline{F(\C^n)} = \C^n$, then $S_F$ is either an empty set or a hypersurface. 
\end{theorem} 

%Note that by $\overline{F(\C^n)}$, we mean the closure of $F(\C^n)$ in $\C^n$. 

\section{The variety $V_F$}  \label{constructionVF}
We recall in this section the construction of the variety $V_F$ and the results obtained in \cite{Valette} and \cite{ThuyValette}: Let $F: \C^n \to \C^n$ be a polynomial mapping. We consider $F$ as a real mapping $F: \R^{2n} \to \R^{2n}$. 
 By $Sing F$ we mean the set of critical points of $F$. 
Thanks to the lemma 2.1 of \cite{Valette}, 
 there exists a covering $\{ U_1, \ldots , U_p \}$ of $M_F = \R^{2n} \setminus Sing(F)$ by semi-algebraic open subsets (in $\R^{2n}$) such that on every element of this covering, the mapping $F$ induces a diffeomorphism onto its image. We may find some semi-algebraic closed  subsets $V_i \subset U_i$ (in $M_F$) which cover $M_F$ as well. By the Mostowski's Separation Lemma (see \cite{Mos}, page 246), for each $ i =1, \ldots , p$, there exists a Nash function $\psi_i : M_F \to \R$,  such that  $\psi_i$ is positive on $V_i$ and negative on $M_F \setminus U_i$. 
 We can choose the Nash functions $\psi_i$ such that $\psi_i (x_k)$ tends to  zero where $\{x_k\}$  is a sequence in $M_F$ tending to infinity. 
We define 
$$V_F : = \overline{(F, \psi_1, \ldots, \psi_p)(M_F)},$$
that means, $V_F$ is the closure of the image of $M_F$ by $(F, \psi_1, \ldots, \psi_p)$. 

The variety $V_F$ is a real algebraic  singular variety of dimension $2n$ in $\R^{2n + p}$, the singular set of which is contained in  $(  K_0(F) \cup S_F ) \times \{0_\R\}^p$, where $K_0(F)$ is the set of critical values  and $S_F$ is the asymptotic set of $F$.  
%We have 
%begin{proposition} [\cite{Valette}]  \label{valette2} 
%{\it Let $F: \C^n \to \C^n$ be a polynomial mapping. There exists a real semi-algebraic pseudomanifold $N_F \subset \R^{\nu}$,
% for some $\nu  = 2n + p$, where $p>0$  such that 
%$$ Sing (N_F) \subset (S_F \cup K_0(F)) \times \{0_{\R^{p}}\},$$
%and there exists a semi-algebraic bi-Lipschitz mapping
%$$h_F: M_F \to N_F \setminus (S_F \cup K_0(F)) \times \{0_{\R^{p}}\}$$
%where $N_F \setminus (S_F \cup K_0(F)) \times \{\{0\}_{\R^{p}}\}$ is equipped with the Riemannian metric induced by $\R^{\nu}.$}
%\end{proposition}

With the construction of $V_F$, we have the two following theorems. Here with the notations $H_*(V)$ ({\it resp.} $IH_*^{\overline{p}}(V)$), we mean the homology ({\it resp.}, the intersection homology) with both compact supports and closed supports.

\begin{theorem} [\cite{Valette}] \label{valette4}
{\it Let $F: \C^2 \to \C^2$ be a polynomial mapping with nowhere vanishing Jacobian. The following conditions are equivalent:

(1) $F$ is non proper,

(2) $H_2(V_F) \neq 0$,

(3) $IH_2^{\overline{p}}(V_F) \neq 0$ for any perversity $\overline{p}$,

(4) $IH_2^{\overline{p}}(V_F) \neq 0$ for some perversity $\overline{p}$.}

\end{theorem}

Form here, we denote by $\hat {F_i}$ the homogeneous component of $F_i$ of highest degree, or {\it the leading form} of $F_i$. 

%We have the following theorem.

\begin{theorem} \cite{ThuyValette} \label{valettethuy2}
{\it Let $F : \C^n \rightarrow \C^n$ be a polynomial mapping with nowhere vanishing Jacobian. 
If {$\rang_{\C} {(D \hat {F_i})}_{i=1,\ldots,n} \geq n-1$}, where $\hat {F_i}$ is the leading form of $F_i$,  then the following conditions are equivalent:

\begin{enumerate}
\item[(1)] $F$ is non proper,
\item[(2)] $H_2(V_F) \neq 0,$
\item[(3)] $IH_2^{\overline{p}} (V_F) \neq 0$ for any (or some) perversity $\overline{p},$
\item[(4)] $IH_{2n-2}^{\overline{p}} (V_F) \neq 0$, for any (or some) perversity $\overline{p}.$
\end{enumerate}
}
\end{theorem}

Notice that the condition ``$F$ is nowhere vanishing Jacobian'' in these theorems means that the set of critical values $K_0(F)$ of $F$ is an emptyset. 

\medskip 

\medskip

The purpose of this paper is to prove that if $F: \C^n \to \C^n$ ($\, n \geq 2$) is a non-proper {\it  generic dominant } polynomial mapping, then the 2-dimensional homology and  intersection homology (with any perversity)  of $V_F$ are not trivial. In order to compute the intersection homology of the variety $V_F$ in the case $K_0(F) \neq \emptyset$,   we have to stratify the set $K_0(F) \cup S_F$. Furthermore, the intersection homology of the variety $V_F$ does not depend on the stratification of $V_F$ if we use a locally topologically trivial stratification. By theorem \ref{doesnotdepend}, a Thom-Mather stratification is a locally topologically trivial stratification. 
In the following section, we provide an explicit 
%\footnote{La m\^eme question que dans la r\'esum\'e.} a 
Thom-Mather stratification of the set $K_0(F) \cup S_F$.

\section{A Thom-Mather Stratification of the set $K_0(F) \cup S_F$}
 
%For the intersection homology of the set $V_F$ does not depend on the stratification of $V_F$, we need a Thom-Mather stratification of the set $K_0(F) \cup S_F$. One way to stratify the set $K_0(F) \cup S_F$ is to stratify $S_F$ and $K_0(F)$. 
% In \cite{Thuy}, by {\it ``la m\'ethode des fa{\c c}ons''} we have a Thom-Mather stratification of the set $S_F$. The following obsevation shows that we can stratify the set $S_F$, but we can not stratify the set $K_0(F)$ alone, so we have to stratify the set $S_F \cup K_0(F)$.

%To compute the intersection homology of the set $V_F$ in the case $K_0(F) \neq \emptyset$,   we have to stratify the set $K_0(F) \cup S_F$. Furthermore, for the intersection homology of the set $V_F$ does not depend on the stratification of $V_F$, we need a locally topologically trivial stratification of the set $V_F$. We know that a Thom-Mather stratification implies a locally topologically trivial stratification (see \cite{Thom}, \cite{Whitney}). 

It is well-known that there exists a partition of the set $K_0(F)$ defined  by Thom  \cite{Thom}. 
% \footnote{la m\^eme question qu'avant sur la r\'ef\'erence}. 
In \cite{Thuy}, we provide a  Thom-Mather stratification of the asymptotic set $S_F$ by  the {\it ``m\'ethode des fa{\c c}ons''}. In this section, we prove that in general the set $K_0(F)$ is not closed, so we cannot define a stratification of $K_0(F)$ satisfying the frontier condition. However, we prove that the set $K_0(F) \cup S_F$ is closed and $K_0(F) \cup S_F = \overline{K_0(F)} \cup S_F$. 
 The fact ``$K_0(F) \cup S_F$ is closed'' has been proved in  \cite{Thuy1}, 
and this result has been generalized in \cite{CidinhaLuis}. 
 Morever, the affirmation  $K_0(F) \cup S_F = \overline{K_0(F)} \cup S_F$ allows us to show that there exists a Thom-Mather stratification of the set $K_0(F) \cup S_F$ compatible with the partition of the set $K_0(F)$ defined by Thom and compatible with the stratification of the set $S_F$ defined by   the {\it ``m\'ethode des fa{\c c}ons''}.

We begin this section by giving an example to show that in general the set $K_0(F)$ of a polynomial mapping $F: \C^n \to \C^n$ is neither closed, nor smooth, nor pure dimensional. Recall that a set $X$ is pure dimensional of dimension $m$ if  any point of this set admits a $m$-dimensional neighbourhood in $X$. 

\begin{example} \label{pasferme}
{\rm 
Let us consider the polynomial mapping 
$F: \C^3_{(x_1,x_2,x_3)} \to \C^3_{(\alpha_1, \alpha_2, \alpha_3)}$ such that 
$$F(x_1,x_2,x_3)= (x_1^3-x_1x_2x_3,x_2x_3,x_3x_1).$$
%Let us denote $ J_F$ by the Jacobian matrix of $F$, we have 
%$$ J_F ={ \begin{pmatrix} 3x_1^2 - x_2x_3 & -x_1x_3 & -x_1x_2  \\ 0 & x_3 & x_2 \\ x_3 & 0 &  x_1
%\end{pmatrix}}.$$ 
Then, the jacobian determinant $\vert J_F(x)\vert$ of $F$  is given by 
$x_1x_3 (3x_1^2 - x_2x_3).$  If $\vert J_F(x) \vert = 0$ then $x_1 = 0$ or $x_3 = 0$ or $3x_1^2 = x_2x_3$.  So we have the following cases:

+ if  $x_1 = 0$ then $F(0, x_2, x_3) = (0, x_2x_3, 0)$ and the axis $0\alpha_2$ is contained in $ K_0(F)$,

+ if $x_3 = 0$ then $F(x_1, x_2, 0) = (x_1^3, 0, 0)$ and the axis $0\alpha_1$ is contained in  $K_0(F)$,

+ if $3x_1^2 = x_2x_3$ then $F(x_1, x_2, x_3) = (-2x_1^3, 3x_1^2, x_3x_1) = (\alpha_1, \alpha_2, \alpha_3)$. We observe that: if $x_1 = 0$ then $\alpha_1 = \alpha_2 = \alpha_3 = 0$; If $x_1 \neq 0$ then $\alpha_1 \neq 0$ and $\alpha_2 \neq 0.$ Moreover, since $3x_1^2 = x_2x_3$ and $x_1 \neq 0$, then $x_3 \neq 0$, this implies $\alpha_3 \neq 0$.  Furthermore, we have $27 \alpha_1^2 = 4\alpha_2^3$. Let 
%Denote  $(\mathscr{S})$ by the part of the surface of equation  
$$(\mathscr{S}) = \{ (\alpha_1, \alpha_2, \alpha_3) \in \C^3_{(\alpha_1, \alpha_2, \alpha_3)}: 27\alpha_1^2 = 4\alpha_2^3,  \alpha_1 \neq 0, \alpha_2 \neq 0, \alpha_3 \neq 0 \},$$
then 
 $(\mathscr{S}) \cup \{ 0 \}$ is contained in $K_0(F)$. 
%such that $\alpha_1 \neq 0, \alpha_2 \neq 0, \alpha_3 \neq 0$. Then, we have 

So, we have $K_0(F) = (\mathscr{S}) \cup 0\alpha_1 \cup 0\alpha_2$ (see figure \ref{dessinK0F}).
\begin{figure}[h!]\begin{center} 
\includegraphics[scale=0.5]{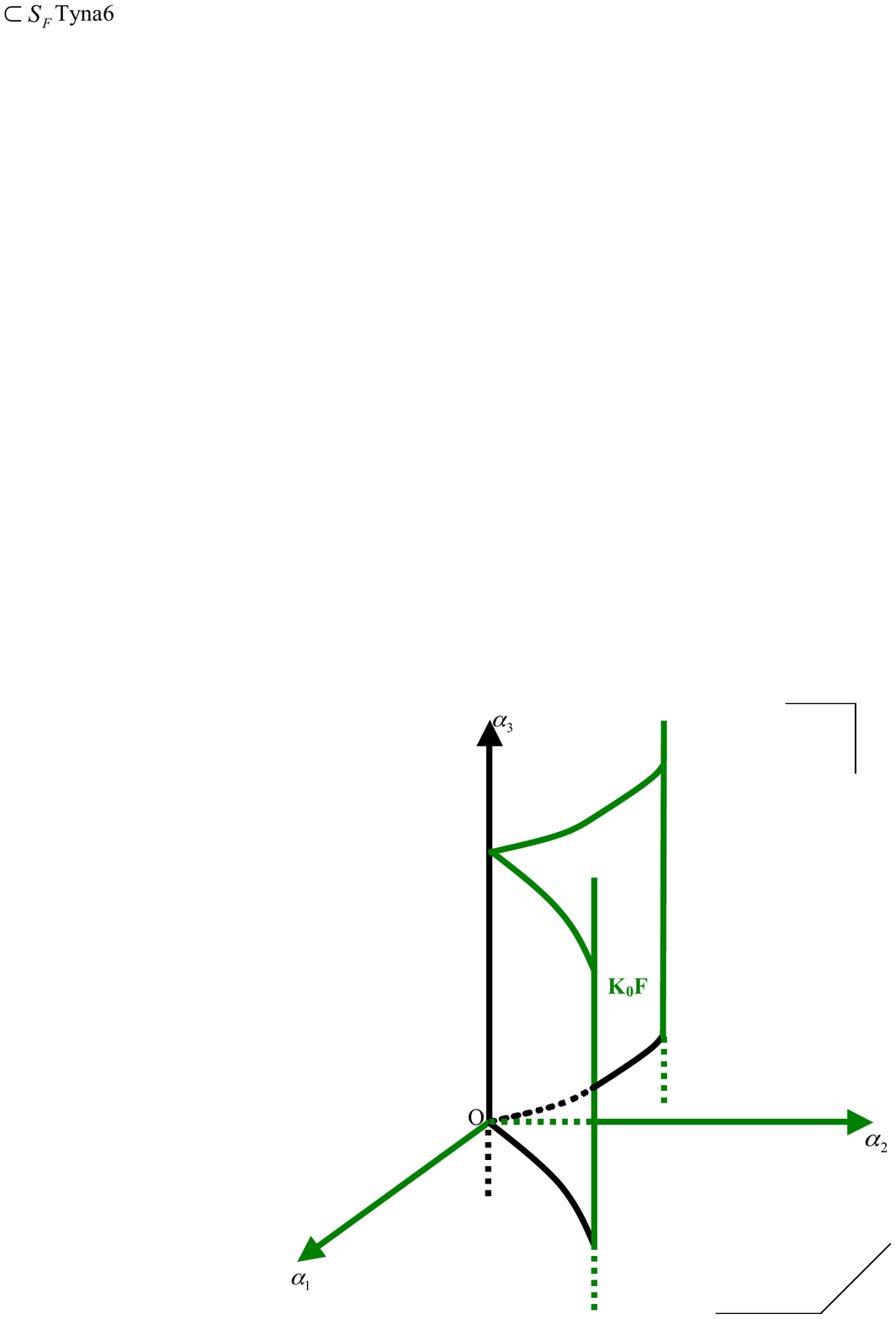}
\caption{The set $K_0(F)$ of the polynomial mapping $F= (x_1^3-x_1x_2x_3,x_2x_3,x_3x_1)$.}
\label{dessinK0F}
\end{center}\end{figure}

Notice that $K_0(F)$ does not contain neither $0\alpha_3 \setminus \{ 0 \}$, nor the curve  $(\mathscr{C})$ of equation  $27\alpha_1^2 = 4\alpha_2^3$ in the plane $(0\alpha_1 \alpha_2)$.
 However $\{0\} \in K_0(F)$ and this is the singular point of $K_0(F)$. So,  the set $K_0(F)$ is  neither closed, nor smooth, nor pure dimensional.
 
}
\end{example}

From the example \ref{pasferme}, in general the set $K_0(F)$ is not closed, so we cannot stratify $K_0(F)$ in such a way that the stratification satisfies  the frontier condition. 
  The following proposition allows us to provide a stratification satisfying  the frontier condition of the set $K_0(F) \cup S_F$.

%L'existence de la stratification de $S_F \cup K_0(F)$ est justifi\'ee par le fait que $S_F \cup K_0(F)$ est ferm\'e~:
\begin{proposition} \label{proK0F}
The set $ K_0(F) \cup S_F $ is closed. Moreover, we have 
$$ K_0(F) \cup S_F  =   \overline{K_0(F)} \cup  S_F. $$
\end{proposition}

To prove this proposition, we need the two following lemmas.
%\begin{lemma} \label{lemJF}
%For a polynomial mapping  $F: \C^n \to \C^n$, the set of the solutions of $\vert J_F (x) \vert = 0$ is closed, where $\vert J_F (x) \vert$ is the jacobian determinant of $F$ at $x$. 
%\end{lemma}
%\begin{proof}
%Considering a sequence  $\{x_k\}$ contained in the set $\{x \in \C^n: \vert J_F(x) \vert = 0 \}$ 
%such that $x_k$ tends to $x_0$. 
%Since $F$ is a polynomial mapping, then $\vert J_F(x) \vert$ is also a polynomial mapping and $\vert J_F (x) \vert $ is continuous. 
%Hence $\vert J_F(x_k) \vert$ tends to $\vert J_F(x_0) \vert $. 
%Since $\vert J_F(x_k) \vert = 0$ for all $x_k$, we have $\vert J_F(x_0) \vert = 0$. 
%So $x_0$ belongs to the set $\{x: \vert J_F(x) \vert = 0 \} $. 
% We conclude that the set of the solutions of $\vert J_F (x) \vert = 0$ is closed.
%\end{proof}

\begin{lemma} \label{lemK0F}
The set $\overline{K_0(F)} \setminus K_0(F)$ is contained in the set $S_F$.
\end{lemma}
\begin{preuve}
Let $a \in \overline{K_0(F)} \setminus K_0(F)$. 
 There exists a sequence $\{a_k\}  \subset K_0(F)$ such that 
 $a_k$ tends to $a$. 
Then there exists a sequence 
 $\{x_k \}$ contained in the set $\{ x: \vert J_F (x) \vert = 0 \}$ such that $F(x_k) = a_k$, for all $k$, where $\vert J_F (x) \vert $ is the determinant of the Jacobian matrix of $F$. 
  Assume that the sequence $\{x_k\}$ tends to $x_0$
 and $x_0$ is finite. 
 Since the set $\{ x \in \C^n: \vert J_F (x) \vert = 0 \} $ is closed, 
  then 
 $x_0$ belongs to the set $\{ x \in \C^n: \vert J_F (x) \vert = 0\}$. 
Moreover, since $F$ is a polynomial mapping, then $F(x_k)$ tends to $F(x_0)$.
 Hence $a_k$ tends to $F(x_0)$ 
and $a = F(x_0)$.  
 Since $x_0$ is finite, then $a \in K_0(F)$, which provides the contradiction. Then $x_k$ tends to infinity and $a$ belongs to $S_F$. 
\end{preuve}

 Considering now the  graph of $F$ in $\C^n \times \C^n$, that means
$$ \graph F=\{(a,F(a)) : a \in \C^n\} \subset \C^n \times \C^n.$$
Let $\overline{\graph F}$ be the projective closure of  $\graph F$ in $\C\P^n \times \C^n$. 
% In order to caracterize the set $S_F$, 
%we consider the closure of $ \graph F$ in $\C\P^n \times \C^n$. 
%Denote $\pi_2$ by the canonical projection $ \C\P^n \times \C^n \to \C^n$,
 We have the following lemma:

\begin{lemma} \label{obverseSF} 
The asymptotic set $S_F$ of a polynomial mapping 
 $F: \C^n \to \C^n$ is the image of the set $\overline{\graph F} \setminus \graph F$ by the canonical projection 
$\pi_2: \C\P^n \times \C^n \to \C^n$.  
\end{lemma}

This lemma is well-known. In fact, this is the first observation of Jelonek \cite{Jelonek} when he studied the geometry of the asymtotic set $S_F$. We can find this fact, for example, in the introduction of \cite{Valette}. We provide here a demonstration of this observation.

\begin{preuve} 
Firstly, we show the inclusion $S_F \subset \pi _2 ( \overline{\graph F} \, \setminus \, \graph F).$ 
Let  $a' \in S_F$, 
there exists a sequence  $\{\xi_k\} \subset \C^n$ such that  
$\xi_k$ tends to infinity and $F(\xi_k)$ tends to $a'.$  
The limit of the sequence $\{(\xi_k, F(\xi_k))\}$ is $a^* = (\infty, a')$, where $a^* \in \overline {\graph F} \setminus \graph F \subset (\C\P^n \times \C^n)$ and $a' =\pi_2(a^*)  \in \pi _2 ( \overline{\graph F} \setminus \graph F).$ 

Now we show the inclusion $\pi _2 ( \overline{\graph F} \, \setminus \, \graph F) \subset S_F.$ 
Let  $a' \in \pi _2 ( \overline{\graph F} \, \setminus \graph F)$,  then 
there exists $a^*=(a, a') \in \overline{\graph F} \setminus \graph F$ 
such that $a^* \in \overline{\graph F}$ but $a^* \notin \graph F$. 
Then we have $a' \neq F(a)$. 
Moreover, there exists a sequence $\{(\xi_k, F(\xi_k))\} \subset \graph F$ 
 such that $(\xi_k,F(\xi_k))$ tends to $(a, a')$. 
Hence the sequence $\xi_k$ tends to $a$ and $F(\xi_k)$ tends to $a'.$ 
Since $F$ is a polynomial mapping, then $F(\xi_k)$ tends to $F(a)$. 
But $a' \neq F(a)$, then $a = \infty$, 
and $\xi_k$ tends to infinity. Thus we have $a' \in S_F.$
\end{preuve}

We prove now the proposition \ref{proK0F}.

\begin{preuve}
By the lemma \ref{obverseSF}, 
the set $S_F$  is the image of the set 
$(\overline{\graph F} \, \setminus \graph F)$ by the canonical  projection $ \pi_2 : \P \C^n \times \C^n \to \C^n.$ Then the set  $S_F$ is closed. Moreover, we have   
$$ \overline{K_0(F)} \cup S_F = K_0(F) \cup (\overline{K_0(F)} \setminus K_0(F)) \cup S_F. $$
By the lemma  \ref{lemK0F}, we have $\overline{K_0(F)} \setminus K_0(F) \subset S_F$, then  
$  K_0(F) \cup S_F  = \overline{K_0(F)} \cup S_F. $
Consequently, the set  $K_0(F) \cup  S_F$ is closed.
\end{preuve}
% \begin{remark}
%{\rm Tous les resultats ci-dessus, dans ce chapitre, sont vrais pour le cas g\'en\'eral $\C^m \to \C^n$.}
%\end{remark}
%En g\'en\'eral, l'ensemble $K_0F$ n'est pas ferm\'{e} et $\overline{K_0F} \setminus K_0F \neq \emptyset$, de plus $(\overline{K_0F} \setminus K_0F) \subset S_F$. La raison pour laquelle  l'ensemble $K_0F$ n'est pas ferm\'{e} est qu'il y a des composantes dans $\overline{K_0F}$ qui ne sont pas contenues dans $F(X)$, alors que $K_0F \subset F(X)$. Ces composantes sont dans $\overline{K_0F}$ mais ne sont pas dans $K_0F$. 
%\footnote{\textcolor{green}{Ici, je voudrais dire ce paragraphe pour insister encore laquelle raison pour que $K_0F$ ne soit pas ferm\'e (mais elle va \^etre ouvert quand  jusqu'\`a il y a une condition ? Tu le sais ? Je sais, parce que j'ai d\'ej\`a essay\'e, tu veux je te montrer ?  C'est incroyable la condition que j'ai trouv\'e. Je vais te montrer avec plaisir.). Mais avec les commentaires que tu m'as propos\'e, il faut l'enlever ce pharagraphe ? - oui, j'ai h\^ate d'avoir Vu, mais tu vas me montrer, n'est-ce pas ?}}

%In the following proposition, we denote by $(V, \Sigma_V)$ a stratification of the variety $V$, 
%where $\Sigma_V$ is the set of all the strata of $V$. 
%Two stratifications $(V, \Sigma_V)$ and $(V', \Sigma_{V'})$ 
%are tranversal if and only if for any stratum $S \in \Sigma_V$ and for any stratum 
%$S' \in \Sigma_{V'}$, then $S$ and $S'$ are tranversal.

\begin{theorem} \label{prostraK0FSF} 
Let $(\mathcal{K})$ be the partition of ${{K_0(F)}}$ defined by Thom in \cite{Thom} and the stra\-ti\-fi\-ca\-tion $(\mathcal{S})$ of $S_F$ defined by  the {\it ``m\'ethode des fa{\c c}ons''} in \cite{Thuy}.  
 Assume that if  $K$ is a stratum of $K_0(F)$ and if $S $ is a stratum  of $S_F$, then  $\overline{K}$  and $S$ are transverse. 
Then there exists a Thom-Mather stratification of the set  $ K_0(F)  \cup S_F$ compatible 
%with the partition of the set $K_0(F)$ defined by Thom and compatible
 with $(\mathcal{S})$ and $(\mathcal{K})$. 
%the stratification of the set $S_F$ defined by   the {\it ``m\'ethode des fa{\c c}ons''}. 
\end{theorem}
\begin{preuve}
By the  proposition \ref{proK0F}, the set  $ K_0(F)   \cup S_F$ is closed and  
$$K_0(F) \cup  S_F  = \overline{K_0(F)} \cup  S_F.$$
% In order to stratify the set 
% $K_0(F) \cup  S_F$, 
% we have to stratify the set $\overline{K_0(F)}  \cap S_F$. 
Con\-si\-de\-ring the partition $(\mathcal{K})$ of ${{K_0(F)}}$ defined by Thom \cite{Thom} and the stratification $(\mathcal{S})$ of $S_F$ defined by  the {\it ``m\'ethode des fa{\c c}ons''} \cite{Thuy}.  
Since  $\overline{K}$ and $S$ are transverse  for any $K \in (\mathcal{K})$ and for any $S \in (\mathcal{S})$, then we can stratify $\overline{K_0(F)} \cap S_F $ by strata $\overline{K} \cap S$ where $K \in (\mathcal{K})$ and $S \in (\mathcal{S})$.
%$$\{ S \cap \overline{K}: S \in (\mathcal{S}), K \in (\mathcal{K}) \}.$$
So we obtain a stratification of the set $ K_0(F) \cup S_F$ compatible 
 with  $(\mathcal{K})$ and $(\mathcal{S})$.  
We have  $\{\overline{K}: K \in (\mathcal{K})\}$ is a Thom-Mather stratification of $\overline{K_0(F)}$ compatible with $(\mathcal{K})$  (\cite{Thom}, theorem 4.B.1) and $(\mathcal{S})$ is a Thom-Mather stratification (\cite{Thuy}, theorem 4.6). 
%\footnote{POURQUOI ??? REFERENCE???}. 
Hence,  
we get a Thom-Mather stratification of the set  $S_F \cup K_0(F)$ (\cite{Tr}, page 4). 

% Assume that 
%A stratification of $S_F \cap \overline{K_0(F)}$ is given by:
%$$\Sigma_{S_F} \cap \Sigma_{\overline{K_0(F)}} = \{ X \cap X' : X \in \Sigma_{S_F}, X' \in \Sigma_{\overline{K_0(F)}}\}.$$
% if $\Sigma_{S_F}$ and $\Sigma_{\overline{K_0(F)}}$ are transversal (see \cite{Tr}, page 4), which provides us the result. 
%\footnote{Verify. Verify if this provides a Thom-Mather stratification of $S_F \cup K_0(F)$.}   
\end{preuve}

%\begin{corollary}
%If $F : \C^n_{x} \to \C^n_{\alpha}$ is dominant, then there exists a natural stratification of  $\C^n_{\alpha}$ by the smooth varieties compatible with la stratification de $S_F$ ((\textcolor{red}{by la m\'ethode des fa{\c c}ons ... ?})) and la partition de $K_0(F)$ (\textcolor{red}{in the proposition ... ?}).  
%\end{corollary}
%\begin{proof}
%Since $F$ is dominante, then we have $\overline{F(\C^n_{x})} = \C^n_{\alpha}$. We see also  $\overline{F(\C^n_{x})} = S_F \cup F(\C^n_{x})$. By the propostion \ref{prostraK0FSF}, we have a stratification of $S_F \cup K_0(F)$ compatible with the partitions of  $K_0(F)$ and $S_F$, a  stratification of $\overline{F(\C_x^n)}$ compatible with $S_F$ and $\overline{K_0(F)}.$ 
%\end{proof}

\medskip 

In order to illustrate the theorem \ref{prostraK0FSF}, let us provide an explicit example. 

\begin{example} 
%\footnote{Modifier cet exemple}
{\rm  Let us consider the example  
\ref{pasferme}: let $F: \C^3_{(x_1,x_2,x_3)} \to \C^3_{(\alpha_1, \alpha_2, \alpha_3)}$ be  the  polynomial mapping   such that  $F(x_1,x_2,x_3)= (x_1^3-x_1x_2x_3,x_2x_3,x_3x_1).$ 

\medskip 

a) At first, via this example, we make clear the idea ``define a partition of the set $K_0(F)$ by {\it constant rank}''   by providing the partition of the set $K_0(F)$ defined by Thom in \cite{Thom}, consisted in the five following steps. 

1) Step 1: Subdividing the singular set $Sing F$ of $F$ into subvarieties $V^i$, where $V^i = \{ (x_1,x_2,x_3) \in \C^3: \Rang J_F(x_1,x_2,x_3) = i\}$. From the example \ref{pasferme}, we have:  
$$\begin{aligned}
V^0 & = \{(0,0,0)\},  \quad V^1 = \{(0,x_2,0) : x_2 \neq 0 \},\cr
V^2 & = \{(0,x_2,x_3) : x_3 \neq 0 \} \cup \{ (x_1,x_2,0) : x_1 \neq 0\} \cup \{ (x_1,x_2,x_3): x_3 \neq 0, x_2x_3 = 3x_1^2\}. \end{aligned}$$
%We know that 
%$$Sing F = \{ (0,x_2,x_3)\} \cup \{ (x_1,x_2,0)\} \cup \{ (x_1,x_2,x_3) : 3x_1^2 = x_2x_3\}.$$
% Let us consider the subvarieties  $V^i = \{ (x_1,x_2,x_3) \in \C^3: \Rang J_F(x_1,x_2,x_3) = i\}$, for $i=0, 1, 2, 3$, we have $SingF = V^0 \cup V^1 \cup V^2$, where 
%$$\begin{aligned}
%V^0 & = \{(0,0,0)\},  \quad V^1 = \{(0,x_2,0) : x_2 \neq 0 \},\cr
%V^2 & = \{(0,x_2,x_3) : x_3 \neq 0 \} \cup \{ (x_1,x_2,0) : x_1 \neq 0\} \cup \{ (x_1,x_2,x_3): x_3 \neq 0, x_2x_3 = 3x_1^2\}. \end{aligned}$$

2) Step 2: Subdividing the sets $V^i$ in step 1 into smooth varieties. Since $V^2$ is not smooth,  so we need to subdivide $V_2$ into  $V_1^2 : = \{ (x_1,x_2,x_3) : 3x_1^2 = x_2x_3, x_3 \neq 0 \}$, $V_2^2 := \{ (0,x_2,x_3) : x_3 \neq 0\}$ and $V_3^2  := \{ (x_1,x_2,0) : x_1 \neq 0\}$.

\medskip 

3) Step 3: Making a partition of the set $Sing F$ from the subsets $V_j^i$ in the steps 1 and 2. 
Since $V_1^2 \cap V_2^2 = 0x_3 \setminus \{ 0 \}$, so  let us consider:
$$\begin{aligned} 
{V'}_1^2  & : = V_1^2 \setminus 0x_3, \quad {V'}_2^2 := V_2^2 \setminus 0x_3, \quad {V'}_3^2 := V_3^2 \setminus 0x_3, \cr
{V'}_1^1 & : = V^1 \setminus \{ 0 \} = 0x_2 \setminus \{ 0 \}, \quad {V'}_2^1 = 0x_3 \setminus \{ 0 \}, \quad
{V'}^0 := \{ 0 \}. \end{aligned}$$
We get a partition of $Sing F$.
%We obtain a partition of $\C^3_{(x_1,x_2,x_3)}$ compatible with $Sing(F)$.
% This partition defines also a stratification of  $\C^3_{(x_1,x_2,x_3)}$
%and $Sing(F)$ (see figure \ref{K0F3Joli}).
%
%\begin{figure}[h!]\begin{center}
%\includegraphics[scale=0.5]{K0F3Joli.pdf}
%\caption{A partition of $\C^3_{(x_1,x_2,x_3)}$ and $SingF$ for the polynomial mapping $F(x_1,x_2,x_3)= (x_1^3-x_1x_2x_3,x_2x_3,x_3x_1)$}
%\label{K0F3Joli}
%\end{center}\end{figure}
%\medskip

\medskip

4) Step 4: Computing $\Rang {J_F}_ {\vert T{V'}_i^j}$. We have
$$ \begin{aligned}
\Rang {J_F}_ {\vert T{V'}_1^2} & = 2, \quad 
 \Rang {J_F}_ {\vert T{V'}_2^2} = \Rang {J_F}_ {\vert T{V'}_3^2} = 1,  \cr
\Rang {J_F}_ {\vert T{V'}_1^1} & = \Rang {J_F}_ {\vert T{V'}_2^1} = \Rang {J_F}_ {\vert T{V'}^0}=0.\end{aligned}$$

5) Step 5: Computing $W^{i,k}_j := F( \{ x  \in {V'}^i_j : \Rang {J_F}_{\vert T_x{V'}^i_j} = k\}).$ 
We have
$$
 W^{2,2} = (\mathscr{S}), \quad 
W^{2,1}_2 = 0 \alpha_2, \quad 
W^{2,1}_3 = 0 \alpha_1, \quad W^{1, 0}_1 = W^{1, 0}_2 = W^{0, 0} = \{ 0 \}.$$
Recall that $(\mathscr{S}) = \{ (\alpha_1, \alpha_2, \alpha_3) \in \C^3_{(\alpha_1, \alpha_2, \alpha_3)}: 27\alpha_1^2 = 4\alpha_2^3,  \alpha_1 \neq 0, \alpha_2 \neq 0, \alpha_3 \neq 0 \}.$ 

\noindent Each $ W^{i,k}_j$ is a $k$-dimensional smooth variety of $K_0(F)$. So we get a partition of $K_0(F)$ by smooth varieties (see figure  \ref{partitionK0F}). 
%Let 
%$${W'}^{i,k}_j : = \{ x \in W^{i,k}_j : x \notin W^{i',k'}_{j'}, \forall (i',j',k') \neq (i, j, k)\}$$
%then $\{ {W'}^{i,k}_j\}$ define a stratification of  $\C^3_{(\alpha_1,\alpha_2,\alpha_3)}$.

\begin{figure}[h!]\begin{center} 
\includegraphics[scale=0.8]{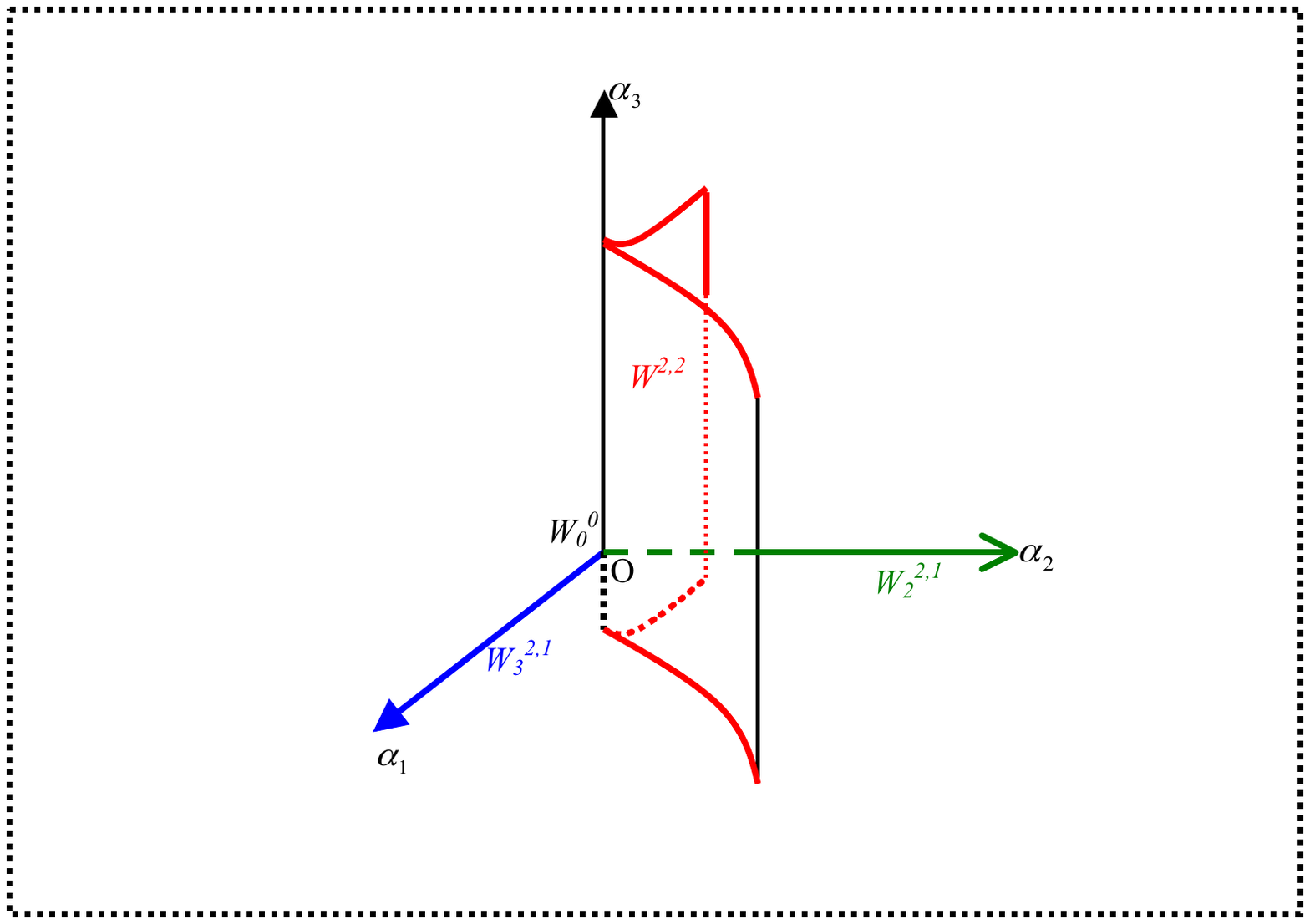}
\caption{The partition of $K_0(F)$ defined by Thom  of the polynomial mapping $$F= (x_1^3-x_1x_2x_3,x_2x_3,x_3x_1)$$}
\label{partitionK0F}
\end{center}\end{figure}

%\footnote{1) Simplifier ce dessin: peut-\^etre ce n'est pas la peine de montrer $K_0(F)$; 2) Change the axis $\alpha, \beta, \gamma$ by $\alpha_1, \alpha_2, \alpha_3$; Chinh lai hinh, luu y phan thay va phan khong thay.}

\medskip 

b) Let us determine now a stratification of $ K_0(F) \cup S_F$. By definition, we can check easily that the asymptotic set $S_F$ is the union of two planes $\{ \alpha_1= 0 \}$ and $\{ \alpha_3= 0 \}$. A  stratification of the set  $ K_0(F) \cup S_F$, which is in fact a Thom-Mather stratification, is given by the filtration:
$$K_0(F) \cup S_F \supset {0\alpha_1 \cup 0\alpha_2}  \supset \{ 0 \} \supset \emptyset,$$
where 

+ $\big( K_0(F) \cup S_F \big) \setminus \ (0\alpha_1 \cup 0\alpha_2) $ has three strata: the two strata $S^2_1 = (0\alpha_2 \alpha_3) \setminus 0\alpha_2$ and $S^2_2 = (0\alpha_1 \alpha_2) \setminus 0 \alpha_2$ defined by {\it fa{\c c}ons} \cite{Thuy} 
and the stratum $ S^2_3 = (\mathscr{S}) \setminus S_F \subset K_0(F) \setminus S_F$ defined by Thom, 

+ $({0\alpha_1 \cup 0\alpha_2})  \setminus \{ 0 \}$ has two strata: the stratum $0 \alpha_1  \setminus \{ 0 \}$ is defined by Thom, the stratum $0\alpha_2  \setminus \{ 0 \}$ is defined either by Thom or by  {\it fa{\c c}ons},

+ \{0\} is the 0-dimensional stratum defined by Thom.

%, where $(\mathscr{C})$ is the curve of equation  $27\alpha_1^2 = 4\alpha_2^3 \subset (0\alpha_1 \alpha_2)$ and $(\mathscr{S})$ is the part of the surface of equations  $27\alpha_1^2 = 4\alpha_2^3, \qquad \alpha_3= t \in \C $
%such that $\alpha_1 \neq 0, \alpha_2 \neq 0, \alpha_3 \neq 0$. 
%The stratum $S^2_3$ is a stratum of  $K_0(F) \setminus S_F$. 

%La strate $0\alpha_2$ a seulement une composante o\`u $\Xi (a) = \{ (2)[1,3] , (3)[1,2] \}$ pour un point $a \in 0\alpha_2.$ 

}

\end{example}

\begin{remark}
{\rm 
%The set $K_0(F)$ is not closed, so we cannot define a stratification of  $K_0(F)$ satisfying the frontier conditions. But the set $S_F \cup K_0(F)$ is closed, so we can define a stratification of the set 
% $S_F \cup K_0(F)$. 
% By {\it ``La m\'ethode des fa{\c c}ons''} in \cite{Thuy}, we have a Thom-Mather stratification of the asymptotic set $S_F$. 
If 
  $K_0(F) \setminus S_F$ is smooth, then we can define easily a stratification of the set  $ K_0(F) \cup S_F$. But in general, $K_0(F) \setminus S_F$ is not smooth. We can check this fact in the following example:
$$F: \C^3 \to \C^3, \quad \quad F= (x_1^3-x_1x_2x_3,x_2x_3,x_3).$$
}
\end{remark}

\begin{remark}
{\rm In all examples in this paper and in \cite{Thuy1}, the set $K_0(F)  \cup S_F$ is pure dimensional if $F$ is dominant. So we can suggest the following conjecture:

\begin{conjecture}
If $F: \C^n \to \C^n$ is a dominant polynomial mapping then the set $K_0(F) \cup S_F$ is pure dimensional.
\end{conjecture} 

Notice that the above conjecture is not true in the case $F$ is not dominant, as shown in the following example: 
$$F: \C^3 \to \C^3, \quad \quad F= (x_1^2 - x_2x_3, x_2 -x_3, x_1-x_3).$$

}
\end{remark}

\pagebreak

\section{The homology and intersection homology of the variety $V_F$  \\ in the case  $K_0(F)$ is non empty}

We have the two following theorems. 

\begin{theorem} \label{Thuy1}
{\it Let $F: \C^2 \to \C^2$ be a non-proper generic dominant polynomial mapping, then 

(1) $H_2(V_F) \neq 0$,

(2) $IH_2^{\overline{p}}(V_F) \neq 0$ for any perversity $\overline{p}$,

(3) $IH_2^{\overline{p}}(V_F) \neq 0$ for some perversity $\overline{p}$.}

\end{theorem}

%Form here, we denote by $\hat {F_i}$ the homogeneous component of $F_i$ of highest degree, or {\it the leading form} of $F_i$. We have the following theorem.

\begin{theorem}  \label{Thuy2}
{\it Let $F : \C^n \rightarrow \C^n$ ($n \geq 3$) be a non-proper generic dominant polynomial mapping. 
If {$\rang_{\C} {(D \hat {F_i})}_{i=1,\ldots,n} \geq n-1$}, where $\hat {F_i}$ is the leading form of $F_i$,  then:

\begin{enumerate}
\item[(1)] $H_2(V_F) \neq 0,$
\item[(2)] $IH_2^{\overline{p}} (V_F) \neq 0$ for any (or some) perversity $\overline{p},$
\item[(3)] $IH_{2n-2}^{\overline{p}} (V_F) \neq 0$, for any (or some) perversity $\overline{p}.$
%\item[(4)] $IH_{2n-2}^{\overline{p}, c} (V_F) \neq 0$, for any (or some) perversity $\overline{p}.$
\end{enumerate}
}
\end{theorem}

%The homology and intersection homology in the theorems \cite{Thuy1} and \cite{Thuy2} are considered with compacted supports or closed supports. 

Before proving these theorems, we recall some necessary definitions and lemmas.

\begin{definition}
A {\it semi-algebraic family of sets (parametrized by $\R$)} 
{\rm is a semi-algebraic set $A \subset \R^n \times \R$, the last variable being considered as parameter.}
\end{definition}
\begin{remark}
{ \rm A semi-algebraic set $A \subset \R^n \times \R$ will be considered as a family parametrized by $t \in \R$. We write $A_t$, for ``the fiber of $A$ at $t$'', {\it i.e.}:
$$ A_t : = \{ x \in \R^n : (x, t) \in A\}. $$}
\end{remark}
\begin{lemma} [\cite{Valette},  lemma 3.1]  \label{valette5} 
{\rm Let $\beta$ be a $j$-cycle and let $A \subset \R^n \times \R$ be a compact semi-algebraic family of sets with $\vert \beta \vert \subset A_t$ for any $t$. Assume that  $\vert \beta \vert$ bounds a $(j+1)$-chain in each $A_t$, $t >0$ small enough. Then $ \beta $ bounds a chain in $A_0$.}
\end{lemma}

\begin{definition}[\cite{Valette}]  \label{coneinfinity}

{\rm Given a subset $X \subset \R^{n}$, we define the {\it ``tangent cone at infinity''}, called ``{\it contour apparent \`a l'infini}''
in \cite{Thuy1} by: 
$$ C_{\infty}(X):=\{\lambda \in \S^{n-1}(0,1) \text{ such that } \exists \eta : (t_0, t_0 + \varepsilon] \rightarrow X \text { semi-algebraic,}$$
$$ \qquad \qquad \qquad \underset{t \rightarrow t_0}{\lim} \eta (t) = \infty, \underset{t \rightarrow t_0}{\lim}\frac{\eta(t)}{\vert \eta(t) \vert} = \lambda \}. $$
}
\end{definition}

\begin{lemma} [\cite{ThuyValette}, lemma 4.10] \label{lemmathuy1}
{\rm Let $F=(F_1, \ldots , F_m) : \R^n \to \R^m$ be a polynomial mapping and $V$ the zero locus of   $\hat{F}: = (\hat{F_1}, \ldots, \hat{F}_m)$, where $\hat{F_i}$ is the leading form of  $F_i$. If $X$ is a subset of $\R^n$ such that $F(X)$ is bounded, then $C_\infty(X)$ is a subset of $\S^{n-1}(0,1) \cap V$, where $V = {\hat{F}}^{-1}(0)$.}
\end{lemma}
%\begin{proof}[Proof]
%By definition, $C_\infty(X)$ is included in $\S^{n-1}(0,1)$. We prove now that  $C_\infty(X)$ is included in $V$. In fact, given $\lambda \in C_\infty(X)$, then there exists a semi-algebraic curve $ \gamma : \, (t_0, t_0 + \varepsilon] \rightarrow X $  such that  $\underset{t \rightarrow t_0}{\lim} \gamma (t) = \infty$ and $\underset{t \rightarrow t_0}{\lim}\frac{\gamma(t)}{\vert \gamma(t) \vert} = \lambda$. Then $\gamma(t)$ can be written as  
%$\gamma(t) = \lambda t^m + \ldots$ and $\hat{G_i} = \hat{G_i}(\lambda)t^{md_i} + \ldots $ where $d_i$ is the homogeneous degree of $\hat{G_i}$.  
%Since $G(X)$ is bounded, then $G_i$ cannot tend to infinity when $t$ tends to $t_0$, hence 
%$ \hat{G_i}(\lambda) =0$  for all $i = 1, \ldots , n$.
%\end{proof}

%Let us prove now Theorem \ref{valettethuy2}. We use the idea and technique 
%of the second and third authors in \cite{Valette}.

\begin{preuve} ({\it Proof of the theorem \ref{Thuy1}}). 
 
The proof of this theorem consists into two steps. In the first step, we use the transversality theorem of Thom (see \cite{Arnold}, page 34): if $F$ is non-proper {\it generic} dominant polynomial mapping,  
we can construct an adapted $(2, \overline{p})$-allowable chain in generic position  providing  non triviality of  homology and  intersection homology of the variety $V_F$, for any perversity $\overline{p}$.   In the second step, we use the same idea than in \cite{Valette} to prove that the chain that we create in the first step cannot bound a  $3$-chain in $V_F$.
\medskip

a) Step 1: Let $F: \C^2 \to \C^2$ be a generic polynomial mapping, then $\dim_{\R} V_F = 4$ (\cite{Valette}, proposition 2.3).  
%the strata of which are the strata of  ${\mathcal{S}}_G \times \{ 0_{\R^{1+p}} \}$ union the strata of the stratification of $F({\mathcal{M}}_G)$ defined by the  rank,  according to Thom (see Remark \ref{remarkstraThom}). 
  Assume that $S_F\neq \emptyset$. 
We claim that $S_F \cap K_0(F) \neq \emptyset$.
%In fact, assume that  $S_F \cap K_0(F) = \emptyset$, that implies $S_F \subset K_0(F)$. 
In fact, since $F$ is dominant, then by the theorem \ref{theoremjelonek1}, we have $\dim_\C S_F = 1$. 
Moreover, since $F$ is generic then $\dim_\C K_0(F) = 1$. Thanks again to the genericity of $F$,  we have $S_F \cap K_0(F) \neq \emptyset$. Let $x_0 \in S_F \setminus K_0(F)$, 
 then there exists a complex Puiseux arc $\gamma$  in $\R^4$, where  
$$\gamma : D(0, \eta) \rightarrow \R^{4}, \quad \gamma = u z^{\alpha} + \ldots, $$
(with $\alpha$ is a negative integer, $u$
 is an unit vector of $\R^{4}$ and $D(0, \eta)$ a small 2-dimensional disc centered in 0 and radius $\eta$)  tending to infinity  in such a  way that $F(\gamma)$ converges to  $x_0$. 
Then, the mapping $h_F \circ \gamma$, where $h_F = (F, \psi_1, \ldots,\psi_p)$ (see the construction of the variety $V_F$, section \ref{constructionVF}) provides a singular $2$-simplex in $V_F$ that we will denote by $c$. We prove now the simplex $c$ is $(\overline{p},2)$-allowable for any perversity $\overline{p}$.  In fact, since  $\dim_\C S_F = 1$, the condition (see \ref{conditionperversity})
\begin{equation*} \label{condition1}
0 = \dim_{\R} \{ x_0 \} = \dim_{\R} (( S_F \times \{ 0_{\R^{p}} \}) \cap \vert c \vert) \leq 2 - \alpha+ p_{\alpha},
\end{equation*}
 where $\alpha = \codim_{\R} S_F = 2$ 
holds for any perversity $\overline{p}$ since  $p_2 = 0$. 

Notice that $V_F \setminus S_F$ is not smooth in general. In fact,  $\Sing (V_F \setminus S_F) \subset K_0(F)$. 
 Let us consider a stratum $V_i$ of the stratification of $K_0(F) \cup S_F$ defined in  the theorem \ref{prostraK0FSF} and denote $\beta = \codim_{\R} V_i.$ 
Assume that $\beta \geq 2$, we can choose the Puiseux arc $\gamma$ such that $c$ lies in the regular part of $V_F \setminus (S_F \times \{0_{\R^{p}}\})$, thanks to the genericity of $F$. In fact, this comes from the generic position of transversality.  
So $c$ is $(\overline{p},2)$-allowable. Hence  we only need to consider the cases $\beta = 0$ and $\beta = 1$. Then:

1) If $c$  intersects $V_i$: since $x_0 \in S_F \setminus K_0(F)$, then we have 
%again by the generic position of transversality, we can choose the Puiseux arc $\gamma$ such that 
$0 \leq \dim_{\R} ( V_i \cap \vert c \vert) \leq 1$. Considering  the condition  
\begin{equationth} \label{condition2}
 \dim_{\R} ( V_i \cap \vert c \vert) \leq 2 - \beta+ p_{ \beta}.
\end{equationth}
We see that $2 - \beta+ p_{ \beta} \geq 1$, for $\beta = 0$  and $\beta = 1$. So the condition (\ref{condition2}) holds. 

2) If $c$  does not meet $V_i$, then the condition 
$$ - \infty = \dim_{\R} \emptyset = \dim_{\R} ( V_i \cap \vert c \vert) \leq 2 - \beta+ p_{\beta}$$
 holds always. 

\noindent In conclusion, the simplex $c$ is $(\overline{p},2)$-allowable for any perversity $\overline{p}$.

% From here, the proof of the theorem follows the ideas of \cite{Valette}. 

We can always choose the Puiseux arc such that the support of $\partial c$ lies in the regular part  of $V_F \setminus (S_F \times \{0_{\R^{p}}\})$ and 
% We have 
%$$H_1(\Reg(V_F \setminus (S_F \times \{0_{\R^{p}}\}))) = 0,$$ 
% then 
 $\partial c $ bounds 
a 2-dimensional singular chain  
$e$ of $\Reg(V_F \setminus (S_F \times \{0_{\R^{p}}\}))$. 
%, where $e$ is a chain with compact supports or closed supports.
So $\sigma = c- e$ is a $(\overline{p},2)$-allowable cycle of  $V_F$. 
%, with compact supports or closed supports.

\medskip 

b) Step 2: We claim that $\sigma$ cannot bound a  $3$-chain in $V_F$.
Assume otherwise, {\it i.e.} assume that there is a  3-chain in $V_F$, satisfying $\partial \tau=\sigma$. Let 
$$A:= h_F^{-1}(\vert \sigma \vert \cap (V_F \setminus (S_F \times \{0_{\R^{p}}\}))),$$ 
$$B:= h_F^{-1}(\vert \tau \vert \cap (V_F \setminus (S_F \times \{0_{\R^{p}}\}))).$$ 
By definition \ref{coneinfinity}, the sets $C_\infty(A)$ and $C_{\infty}(B)$ are subsets of $\S^{3}(0,1)$. 
 Observe that, in a neighborhood of infinity,  $A$ coincides with the support of the  Puiseux arc $\gamma$.  The set $C_\infty(A)$ is equal to $\S^1.a$ 
(denoting the orbit of $a \in \C^2$ under the action of $\S^1$ on $\C^2$, $(e^{i\eta},z) \mapsto  e^{i\eta}z$). 
Let $V$ be the zero locus of the leading forms  $\hat{F}: = (\hat{F_1}, \hat {F_2})$.  
Since $F(A)$ and $F(B)$ are bounded, then by  lemma \ref{lemmathuy1}, the sets $C_\infty(A)$ and $C_{\infty}(B)$ are subsets of $V \cap \S^{3}(0,1).$

 For $R$ large enough, the sphere $\S^{3}(0, R)$ with center 0 and radius $R$ in $\R^{4}$ is transverse to $A$ and $B$ (at regular points). Let 
$$\sigma_R : = \S^{3}(0 , R) \cap A, \qquad \tau_R : = \S^{3}(0, R) \cap B.$$
 Then $\sigma_R$ is a chain bounding the chain $\tau_R$. 
Considering a semi-algebraic strong deformation retraction 
$\Phi : W \times [0;1] \rightarrow \S^1.a$, where $W$ is a neighborhood of  $\S^1.a$ in $\S^{3}(0,1)$ onto $\S^1.a$. Considering $R$ as a parameter, we have the following semi-algebraic families of  chains: 

1) ${\tilde{\sigma }}_R :=\frac{\sigma _R}{R},$
for $R$ large enough, then ${\tilde{\sigma}}_R$ is contained in $W$,

2)  $\sigma'_R = \Phi_1({\tilde{\sigma}}_R)$, where $\Phi_1(x) : = \Phi(x, 1), \qquad x \in W$,

3) $\theta_{R} = \Phi({\tilde{\sigma}}_R)$, we have 
 $\partial \theta_R  = \sigma'_R - {\tilde{\sigma}}_R ,$

4) $\theta'_{R} = \tau_R + \theta_R$, we have 
 $\partial \theta'_R  =  \sigma'_R.$

\noindent As, near infinity, $\sigma_R$ coincides with the intersection of the support of the arc $\gamma$ with $\S^{3}(0, R)$, for $R$ large enough the class of $\sigma'_R$ in $\S^1.a$ is nonzero.

Let $r = 1/R$, consider $r$ as a parameter, and let $\{ {\tilde{\sigma}}_r \}$, $\{ \sigma'_r \}$, $\{ \theta_r \}$ as well as $\{ \theta'_r \}$ the corresponding  semi-algebraic families of  chains. 

Let us denote by $E_r \subset \R^{4} \times \R$ the closure of $|\theta_r|$, and set $E_0 : = (\R^{4} \times \{ 0 \}) \cap E$. Since the strong deformation retraction  $\Phi$ is the identity on $C_{\infty}(A) \times [0,1]$, we see that  
$$E_0 \subset \Phi(C_{\infty}(A) \times [0,1]) = \S^1. a \subset V \cap \S^{3}(0,1).$$

Let us denote by $E'_r \subset \R^{4} \times \R$ the closure of $|\theta'_r|$, and set  $E'_0 : = (\R^{4} \times \{ 0 \}) \cap E'$. Since $A$ bounds $B$, then $C_{\infty}(A)$ is contained in $C_{\infty}(B)$. We have 
$$E'_0 \subset E_0 \cup C_{\infty}(B) \subset V \cap \S^{3}(0,1).$$

The class of $\sigma'_r$ in $\S^1.a$ is, up to a product with a nonzero constant, equal to the generator of $\S^1.a$. Therefore, since $\sigma'_r$ bounds the chain $\theta'_r$, the cycle $\S^1.a$ must bound a chain in $|\theta'_r|$ as well. By Lemma \ref{valette5}, this implies that $\S^1.a$ bounds a chain in $E'_0$ which is included in $V \cap \S^{3}(0,1)$.

 The set $V$ is a projective variety which is an union of cones in $\R^{4}$. Since   $\dim _{\C} V \leq 1$, so $\dim_{\R} V \leq 2 $ and  $\dim_{\R} V \cap \S^{3}(0,1) \leq 1.$
The cycle $\S^1 . a$ thus bounds a chain in $E'_0\subseteq V \cap \S^{3}(0,1)$, which is a finite union of circles, that provides a contradiction.
%are subsets of $V \cap \S^{2n-1}(0,1)$ which are union of circles in  $\R^{2n}$. The set $C_{\infty}(A)$ is a circle in $V \cap \S^{2n-1}(0,1)$, that provides a contradiction. 
\end{preuve}

\begin{preuve} ({\it Proof of the theorem \ref{Thuy2}}). 

Assume that $F : \C^n \rightarrow \C^n$ ($n \geq 3$) is a non-proper generic dominant polynomial mapping. 
Similarly to the previous proof, we have: 

$\bullet$ Since $F$ is dominant, then by the theorem \ref{theoremjelonek1}, we have $\dim_\C S_F = n-1$.  
Moreover, since $F$ is generic then $\dim_\C K_0(F) = n-1$. Thanks again to  the genericity of $F$,  we have $S_F \cap K_0(F) \neq \emptyset$. 
Let $x_0 \in S_F \setminus K_0(F)$, 
 then there exists a complex Puiseux arc $\gamma$  in $\R^{2n}$, where 
$$\gamma : D(0, \eta) \rightarrow \R^{2n}, \quad \gamma = u z^{\alpha} + \ldots, $$
(with $\alpha$ is a negative integer and $u$
 is an unit vector of $\R^{2n}$)  tending to infinity  such a  way that $F(\gamma)$ converges to  $x_0$. Since $x_0 \in S_F \setminus K_0(F)$ and $F$ is generic, then we can choose the arc Puiseux $\gamma$ in generic position, that means the simplex $c$ is $(\overline{p},2)$-allowable for any perversity $\overline{p}$. 

$\bullet$ Now, with the same notations than the above proof, we have:  Since $\Rank_{\C} {(D \hat {F_i})}_{i=1,\ldots,n} \geq n-1$ 
then 
$ \corank_{\C}  {(D \hat {F_i})}_{i=1,\ldots,n}  \leq 1$. Moreover since $\dim _{\C} V = \corank_{\C}  {(D \hat {F_i})}_{i=1,\ldots,n}$ then  $\dim_{\R} V \leq 2 $ and $\dim_{\R} V \cap \S^{2n-1}(0,1) \leq 1.$ The cycle $\S^1 . a$ bounds a chain in $E'_0\subseteq V \cap \S^{2n-1}(0,1)$, which is a finite union of circles, that provides a contradiction. 
%So we have 
%$$IH_2^{\overline{p}, c}(V_F, \R) \ne 0, \quad IH_2^{\overline{p}, cl}(V_F, \R) \ne 0, \quad 
% H_2^{c}(V_F, \R) \ne 0  \text{ and } H_2^{cl}(V_F, \R) \ne 0.$$  

\noindent Hence, we get the facts (1) and (2) of the theorem. Moreover, from the  Goresky-MacPherson Poincar\'e  duality theorem (theorem \ref{ThGM}), we have
$$IH_2^{\overline{p}, c}(V_F) = IH_{2n-4}^{\overline{q},cl}(V_F)$$
where  $\overline{p}$ and $\overline{q}$ are   complementary perversities.  
Since the chain $\sigma$ that we create in the proof of the theorem \ref{Thuy1} can be either a chain with compact supports or a chain with closed supports, so we get the fact (3) of the theorem.

%At first, notice that since $F$ is dominant, then by the theorem \ref{theoremjelonek1}, we have $\dim_\C S_F = n-1$. We claim that $S_F \cap K_0(F) \neq \emptyset$.
%%In fact, assume that  $S_F \cap K_0(F) = \emptyset$, that implies $S_F \subset K_0(F)$. 

%Then the proof of this theorem is similar to the proof of the theorem 4.5 in \cite{ThuyValette}:  

%\text{ and }  IH_2^{\overline{t}, cl}({\mathcal{V}}_G, \R) = IH_{2n-4}^{\overline{0}, c}({\mathcal{V}}_G, \R),$$ 
%that implies 
% $H_{2n-4}^{c}({\mathcal{V}}_G, \R) \ne 0$ and $H_{2n-4}^{cl}({\mathcal{V}}_G, \R) \ne 0$.
%%since 
%% $IH_{2n-4}^{cl, \overline{0}}({\mathcal{V}}_G, \R) = H_{2n-4}^{cl}({\mathcal{V}}_G, \R).$
% 
\end{preuve}

%\textcolor{red}{INTRODUCTION HERE FOR THE NEXT SECTION, EXPLAIN WHY WE NEED A THOM MATHER-STRATIFICATION OF $K_0(F) \cup S_F$ ???}

\begin{remark}
{\rm The properties of the homology and intersection homology in the theorem \ref{Thuy1} and \ref{Thuy2} hold for both compact supports and closed supports.}
\end{remark}

\begin{remark}
{\rm From the proofs of the  theorems \ref{Thuy1} and \ref{Thuy2}, we see that  the properties of the intersection homology in these  theorems do not hold if $F$ is not dominant. The reason is that the theorem \ref{theoremjelonek1} is not true if $F$ is not dominant and then the condition   (\ref{condition2}) does not hold. However, the properties of  the homology hold even if $F$ is not dominant. So we have the two following corollaries. 
}
\end{remark}

\begin{corollary} 
{\rm Let $F: \C^2 \to \C^2$ be a non-proper generic  polynomial mapping, then  $H_2(V_F) \neq 0$.}
\end{corollary}

%Form here, we denote by $\hat {F_i}$ the homogeneous component of $F_i$ of highest degree, or {\it the leading form} of $F_i$. We have the following theorem.

\begin{corollary}  
{\rm Let $F : \C^n \rightarrow \C^n$ ($n \geq 3$) be a non-proper generic polynomial mapping. 
If {$\rang_{\C} {(D \hat {F_i})}_{i=1,\ldots,n} \geq n-1$}, where $\hat {F_i}$ is the leading form of $F_i$,  then $H_2(V_F) \neq 0$.}
\end{corollary}

\begin{remark}
{\rm In the previous papers \cite{Valette} and  \cite{ThuyValette}, the condition ``$F$ is nowhere vanishing Jacobian''  (see theorems \ref{valette4} and \ref{valettethuy2}) implies $F$ is dominant.
 Hence, the condition ``$F$ is dominant'' in the theorems \ref{Thuy1} and \ref{Thuy2}  guarantees  the condition of dimension of the set $S_F$ (see theorem \ref{theoremjelonek1}). Moreover, we need this condition in this paper also to be free ourself from the condition $K_0(F) = \emptyset$, since the condition of dimension of $S_F$ when $F$ is dominant  also guarantees the (generic) tranversal position of  the $(2, \overline{p})$-allowable chain which   provides  non triviality of  homology and  intersection homology of the variety $V_F$ when $K_0(F) \neq \emptyset$ in theorems \ref{Thuy1} and \ref{Thuy2}.  

%
% That is the reason why we need the condition ``$F$ is dominant'' in the theorems \ref{Thuy1} and \ref{Thuy2} to guarantee the condition of dimension of the set $S_F$, thanks to the theorem \ref{theoremjelonek1}. Moreover, the condition of dimension of $S_F$ when $F$ is dominant   guarantee also the (generic) tranversal position of  the $(2, \overline{p})$-allowable chain which   provides  non triviality of  homology and  intersection homology of the variety $V_F$, for any perversity $\overline{p}$ in theorems \ref{Thuy1} and \ref{Thuy2}) when $K_0(F) \neq \emptyset$.  
 }
\end{remark}

\bibliographystyle{plain}

\end{document}